# Sobre o papel dos Departamentos de Matemática na vida e desenvolvimento da comunidade.

Este artigo pretende estimular a discussão no meio matemático sobre o lugar que a Universidade e, em particular, os respectivos Departamentos de Matemática (DM) deve ocupar na comunidade em que se insere. A perspectiva clássica de que a Universidade, por ser portadora do conhecimento científico, desempenha um papel fulcral no desenvolvimento intelectual é hoje questionada pelo simples facto de que a sua participação no processo educativo global é apenas uma pequena parcela de entre muitas componentes.

Criado em 2005 por docentes do departamento de Matemática da Universidade de Aveiro (DMat), o grupo de trabalho EECM [4] dedica-se desde o início a desenvolver actividades na área da educação não formal e que são implementadas no exterior do DMat. Neste processo é envolvido um largo espectro da comunidade: crianças dos 4-13 anos, pais e educadores dessas crianças, jardins-de-infância e escolas de todos os ciclos do ensino básico e do ensino secundário. Convém aqui salientar a organização de "workshops" de matemática destinados a pais/educadores de crianças dos 4-9 anos. Foi também da responsabilidade do EECM a organização, em 2007, dos dois primeiros cursos de matemática da Academia de Verão da UA dirigidos a alunos dos grupos etários 12-15 e 15-18, [4].

Um dos problemas com que nos deparámos desde sempre relaciona-se com o entendimento por parte dos nossos colegas sobre a natureza do trabalho por nós desenvolvido. Muitos são aqueles para os quais a educação não formal é uma actividade menor, e que consideram ser uma "brincadeira" todo o envolvimento com crianças, daí não ser uma actividade importante... A existência deste preconceito e tentar alterar esta situação é a principal razão deste artigo.

Tendo presente a actividade por nós desenvolvida até ao presente momento e que pretendemos continuar a desenvolver, parece-nos primordial esclarecer quais os princípios que norteiam a nossa acção. Mais, esperamos deste modo poder sensibilizar alguns dos nossos pares para também eles desenvolverem trabalho em educação não formal com/para crianças ou com/para alunos dos ensinos básico e secundário. Contrariamente ao que possa parecer à primeira vista, a nossa acção não é motivada por razões quixotescas próximas de uma militância altruísta. O que nos move é acreditar com convicção nos resultados que deste modo é possível alcançar, [8].

Vamos começar por expor o que pensamos serem alguns dos problemas da actualidade que atravessam transversalmente todas as universidades. Simultaneamente, vamos também tentar explicar porque é que a actuação dos DM é importante e porque é que cremos estarem reunidos os pressupostos para que eles desempenhem um papel significativo junto da comunidade.

Um dos aspectos fundamentais para as instituições é o seu financiamento. Actualmente, o financiamento das universidades depende em grande parte do afluxo de novos estudantes. Daí ser fundamental para as instituições cativar os futuros alunos… Deste modo, torna-se importante a opinião que a comunidade tem sobre a capacidade que a universidade demonstra para preparar especialistas competitivos nas respectivas áreas de actividade: com boas bases de formação e com bons hábitos de trabalho.

No panorama actual, o estudante médio que ingressa no ensino superior revela uma preparação desadequada para os programas que lhe são propostos no âmbito da matemática superior, manifestando frequentemente grandes dificuldades perante os mesmos. Contrariando todas as expectativas existentes, a realidade revela que o estudante que chega à universidade foi sujeito, durante a sua escolaridade, a processos que incidem preferencialmente na repetição de técnicas e na execução de sequências de diferentes algoritmos. (O leitor interessado na problemática do behaviourismo em educação encontra bastante literatura disponível, nomeadamente [5]).

Durante décadas um dos objectivos principais do ensino secundário consistia em preparar os alunos para o seu ingresso no ensino superior. Actualmente, entre as finalidades do ensino secundário já não se encontra tal desígnio. Não é de admirar, portanto, que exista um desfasamento entre a preparação dos jovens à entrada nas universidades e a exigência dos programas universitários. Urge resolver este problema.

**O que temos e o que propomos:**

1. O ensino universitário depende da qualidade da preparação dos jovens que concluem o ensino secundário,
2. As escolas secundárias não têm total responsabilidade sobre essa mesma qualidade.

**Conclusão:** À universidade interessa-lhe poder controlar e poder influenciar todo o processo de qualificação dos estudantes.

Não nos parece correcto tentar atribuir a causa do problema da preparação dos jovens que entram no ensino superior à qualidade da formação inicial dos professores ou à quantidade das acções de formação contínua que frequentam durante a sua carreira profissional. Desde o primeiro ciclo que os professores trabalham com programas que foram pensados para turmas com grande número de crianças. O dimensionamento destas turmas não permite que se desenvolva com regularidade trabalho individual com as crianças. Ora, este tipo de trabalho é fundamental na detecção da existência de focos de incompreensão e na eliminação dos mesmos.

A existência de círculos experimentais de matemática (CEM) a funcionar regularmente nas escolas e com a participação dos professores das universidades pode contribuir para a resolução deste problema [1]. A frequência de cursos intensivos durante as férias de verão (ou outros períodos de férias) permite que os jovens aprofundem alguns aspectos do conhecimento que não lhes é proporcionado no ensino regular.

Preparar uma geração é um processo demorado e envolve um investimento que tem que ser feito ao longo de um vasto período. Por vezes, o entendimento entre financiadores e universidades não é o melhor. O apoio económico depende, entre outros factores, da opinião que a sociedade tem da importância do trabalho desenvolvido pelos cientistas.

***O ensino de crianças por cientistas: gasto desproporcional de recursos?***

Começando por um exemplo: há pouco tempo uma mãe colocou a seguinte situação: as suas filhas gémeas, de 6 anos de idade, perguntaram-lhe, "*porque é que os números não param de aumentar? Existe um último número?*" A mãe em causa tem uma boa formação pessoal e receou dar uma resposta que fosse incorrecta. Note-se que esta mesma questão colocada a um cientista pode ser motivo para discutir com as crianças muitos modelos e teorias simples. Podendo o tema ser abordado sucessivamente ao longo dos anos. Pode começar-se com a sugestão de se (re)contarem todos os bebés nascidos e/ou que ainda irão nascer; depois pode pedir-se às crianças que sugiram elas próprias algum assunto. Ao analisar conjuntamente o tema proposto pela criança, pode acontecer que este último conduza a uma situação em que exista um "último número". A regularidade de tais conversas permite que depois de alguns anos ela chegue a discutir, por exemplo, a incomensurabilidade da diagonal do quadrado. Para o cientista (e para o professor com experiência de cientista) cada pergunta colocada pela criança torna-se num tópico para uma conversa, durante a qual ele pode fazer observações ou colocar questões que estejam na base de futuras reflexões. Ainda relativamente ao exemplo colocado, pode perguntar-se*: Quantos milhares são precisos para contar as estrelas? Mil? Pode ser... Porque se acendem as estrelas, ... Se as estrelas se acendem, alguém deve precisar disso…*

A sociedade portuguesa manifesta preocupação e alguma dificuldade em aceitar o analfabetismo e a iliteracia. Surpreendentemente, o mesmo não acontece relativamente ao desconhecimento de conceitos elementares associados à matemática. Não existe qualquer relutância não só em aceitar, como até em desculpabilizar a iliteracia matemática. O mais preocupante é que esta atitude atravessa transversalmente toda a sociedade portuguesa, incluindo os próprios professores. Quanto aos pais, a sua atitude é mais de resignação, vendo este problema como mais um dos desígnios do nosso fado! Em geral, também os próprios pais não tiveram grande sucesso em matemática!

Lutar contra esta situação é uma das finalidades da realização de workshops para pais de crianças em idade do pré-escolar e do 1º ciclo do ensino básico. Durante estas sessões de trabalho é possível (in)formar os pais sobre algumas ideias presentes no processo de educação integral das crianças, desde a mais tenra idade, e explicar-lhes o valor fundamental que a colaboração dos pais tem nesse processo[4\workshop]. A sociedade condena todo o pai que não providencie o bem estar físico dos seus filhos, sendo de aceitação comum que esta é uma das obrigações dos pais. No entanto, a mesma sociedade manifesta-se indiferente para com aqueles pais que deixam os seus filhos sem alimento intelectual. Mudar a atitude dos pais sobre a importância real do papel que podem ter no desenvolvimento intelectual dos seus filhos é uma tarefa fundamental e que urge ser implementada.

***Qual é o objectivo do cientista que trabalha com crianças e de que forma é que ele pode ajudar os professores?***

O sistema de educação actual tem por base o axioma de Volter :" *a criança é uma folha de papel em branco que precisamos de preencher*". Assim, os professores fornecem às crianças vários conjuntos de saberes. Convém salientar, aqui, a importância de cada falha neste sistema. Por exemplo, os danos provocados pela acção de um mau professor durante um certo período de tempo podem levar alguns anos a ser reparados. Se recorrermos às ideias de Pestalozzi e tornarmos interessante o objecto de estudo, o processo educativo torna-se mais efectivo. Quando as crianças têm interesse nos assuntos, elas próprias solicitam os saberes. Ao ajudarmos a criança a sentir em si a capacidade de analisar, raciocinar e concluir, fazemos com que ela deixe de reproduzir acriticamente aquilo que vê. Neste caso já não são os professores que reclamam pelos conhecimentos dos alunos, mas é a própria criança que fica na posição de reclamante. Em termos futuros podemos ter a certeza que, independentemente das circunstâncias, uma criança cuja educação aposte no desenvolvimento dos saberes que a interessam, tornar-se-á num bom especialista. Neste sentido, a actuação do cientista vai permitir despertar e aumentar o interesse da criança pela descoberta e pelo conhecimento. A sua acção não pretende de modo algum substituir o papel do professor em situação de sala de aula, mas sim ajudar a descobrir novos focos de interesse para a criança, abordando os respectivos assuntos de forma clara e com a profundidade adequada ao nível etário e desenvolvimento da criança.

***A dificuldade em trabalhar com alunos do pré-escolar.***

O trabalho de educação não formal com crianças do pré-escolar permite ao professor universitário testar o seu próprio sistema de comunicação. As crianças desta idade são espontâneas, estão abertas para coisas novas e não têm juízos pré formados. A espontaneidade, em particular, permite ter a certeza de que se as crianças não o querem ouvir, isto deve-se ao facto de o assunto que lhes está a apresentar não os interessar. Para interessar a criança é necessário “descer” ao seu nível e é preciso conhecer muito bem as representações que ela tem. Isto só é possível quando se realiza trabalho individual com a criança. “Trabalho” é a palavra adequada, porque a preparação das sessões com crianças exige uma constante reflexão sobre o que cada uma sabe e quais são as suas respectivas representações. É fundamental que se tente perceber porque em dado momento se desinteressam as crianças, isto no sentido de averiguar se a razão do desinteresse reside em alguma falha de comunicação ou se é o assunto que é inadequado. Para que as crianças se interessem pelos assuntos não é necessário transformar tudo em brincadeira e que o seu interlocutor seja um “palhaço”, mas é necessário que elas o entendam.

### *Existe um sentido que determine a interacção dos DM e das escolas?*

Ao longo da história pode observar-se que quando uma sociedade se apercebe da existência de defeitos no sistema educativo procede de imediato à criação de comissões e por vezes à atribuição de fundos. Na realidade, a educação escolar existe já há mais de 2 mil anos e, ingenuamente, pensa-se sempre que mais um esforço vai conduzir o sistema para a perfeição. (Um exemplo clássico é o “boom” das reformas que foram efectuadas nos EUA imediatamente após o lançamento do primeiro satélite russo [5,6]).

O autor (E.L.) na altura em que iniciou a sua actividade de educação complementar em matemática na Rússia ficou surpreendido por não existir uma estratégia comum na forma de lidar com as crianças. Em cada região do país foi desenvolvida uma abordagem própria. Isso foi devido não só à vastidão do país, mas também ao facto de muitos dos cientistas e professores terem motivações diferentes. No entanto, todos eles tinham uma causa comum: o trabalho com crianças. Um outro aspecto muito importante consiste na ajuda dada pelos alunos "mais velhos" e estudantes universitários. Este tipo de experiência está presente na sociedade portuguesa: organizações de escuteiros trabalham recorrendo a um princípio semelhante há já mais de cem anos. O processo de funcionamento é bastante simples: os alunos "mais velhos" e os estudantes universitários participam como monitores dando assistência ao processo de educação dos mais jovens. Deste relacionamento resulta a criação de laços afectivos entre as crianças e os jovens monitores. Estes tornam-se modelos para as crianças e, estas com o decorrer do tempo criam o desejo de ser como eles. Saliente-se que o recurso a tal prática pedagógica é muito útil para a formação de futuros especialistas, sem que para isso seja necessário gastar muito tempo. O recurso a monitores possibilita, também, o envolvimento de um maior número de crianças.

### **Alguns problemas sociais que diminuem o nível intelectual de um país.**

Não podemos deixar de referir alguns aspectos sociais que estão relacionados com esta problemática. Começamos por enumerar alguns problemas sociais que de modo indirecto afectam os interesses da sociedade matemática. Embora não seja objecto do nosso trabalho propor soluções para estes problemas, como cientistas sabemos que a formalização de um problema é já parte da solução.

1. Actualmente é normal colocar os bebés nas creches a partir dos 4 meses de idade. Os pais têm pouco tempo diário para “gastar” com os seus filhos. Esta situação acaba por influenciar a qualidade intelectual dos futuros licenciados das escolas. Vamos aqui referir uma experiência com um grupo de crianças órfãs descrita em [2]. A cada bebé foi atribuída uma “mãe” substituta. As “mães” eram mulheres com graves problemas mentais alheamento da realidade, problemas cognitivos, etc.) que viviam no mesmo orfanato. Tudo

o que elas podiam fazer era abraçar os bebés, mudar-lhes as fraldas, brincar com eles. Não lhes transmitiam nenhum conhecimento. Ora, nesta idade o calor emocional é mais importante do que todas as restantes formas de actividade e educação. Passados 2 anos os indicadores mostraram que o nível intelectual destas crianças aumentou em média cerca de 20-30 pontos, enquanto que, na mesma altura, o nível intelectual das crianças do grupo de controlo foi descendo. O nosso objectivo não consiste em criticar a situação real, apenas constatamos e formalizamos um problema que existe.

2. Os professores das escolas não podem entregar-se completamente à sua tarefa de preparar bem os alunos e não podem ter disponibilidade temporal e emocional para colaborar com os cientistas se as suas necessidades mínimas não estiverem satisfeitas. Não nos referimos às condições salariais, mas sim à estabilidade de trabalho. Só havendo estabilidade é que o professor pode educar bem, e neste caso ele vai tornar-se num modelo para o aluno. A sociedade precisa de lhe garantir condições mínimas para o exercício da sua actividade: possibilidade de viver com a família, de morar perto do local de trabalho e não ter que gastar em viagens 4 horas, ou mais, por dia. Se um professor se esgota em questões exteriores ao ensino e não se realiza na sua profissão, as crianças não o vão respeitar e não vão acreditar no que ele ensina. Além disto, este problema não deixa tempo livre nem disposição ao professor para, por exemplo, colaborar com os cientistas na criação de círculos experimentais.

**Referências:**

[1] S.Genkin, I.Itenberg, D.Fomin, Mathematical Circle: (Russian Experience), American Mathematical Society; 1 edition (July 1996)
[2] Urie bronfenbrenner, Two worlds of childhood: U.S. and U.S.S.R. New York: Russel Sage Foundation 1970.
[3] http://www.ams.org/distribution/mmj/vol1-4-2001/dedication.html
[4] Escola de Educação Complementar em Matemática, http://eecm.mat.ua.pt
[5] Malaty, G., 1998, Eastern and Western mathematical education: unity, diversuity and problems. In: Int. J.Math. Educ/ Sci. Technol., 3, pp/ 421-436
Solmu, 1997.
[6] R.Feinman, Surely You're Joking, Mr. Feynman! W. W. Norton, 1985
[7] www.naturalmath.com, Maria Droujkova.
[8] C.P.Peirce, Issues of pragmatism, Dover Publications, Inc, 1955.